\documentclass{amsart}
\usepackage{amssymb,amsfonts,latexsym}
\usepackage[cspex,bbgreekl]{mathbbol}
\usepackage{bbm}
\usepackage[cp1251]{inputenc}
\usepackage[english,russian]{babel}
\usepackage[matrix,arrow,curve]{xy}
\usepackage{hyperref}
\usepackage{tikz}
 \usetikzlibrary{arrows}
 \usepackage{verbatim}
\newtheorem{theorem}{Теорема}[section]

\newtheorem{proposition}[theorem]{Предложение}
\newtheorem{corollary}[theorem]{Следствие}
\theoremstyle{definition}
\newtheorem{definition}[theorem]{Определение}
\newtheorem{example}[theorem]{Пример}
\theoremstyle{remark}
\newtheorem{remark}[theorem]{Замечание}
\newtheorem{question}[theorem]{Вопрос}
\theoremstyle{remark}

\numberwithin{equation}{theorem}
\renewcommand{\O}{\operatorname{O}}
\newcommand{\Gram}{\operatorname{Gram}}

\newcommand{\rk}{\operatorname{rk}}
\newcommand{\y}{\operatorname{\texttt{ё}}}
\newcommand{\yz}{\operatorname{\texttt{ёж}}}
\newcommand{\Y}{\operatorname{\texttt{Ё}}}
\newcommand{\Pol}{\operatorname{Pol}}
\renewcommand{\S}{\operatorname{S}}

\newcommand{\R}{\mathbbm{R}}

\renewcommand{\Im}{\operatorname{Im}}
\newcommand{\Z}{\mathbbm{Z}}
\newcommand{\N}{\mathbbm{N}}

\newcommand{\eps}{\varepsilon}

\newcounter{itemnumber}

\begin{document}
\sloppy
\author[Н.А.~Печенкин]%
{Николай А. Печенкин}
\title{Ортантные полиэдры}
\address{отдел алгебраической геометрии МИАН}
\email{kolia.pechnik@gmail.com}
\begin{abstract}
В работе изучается возможность реализации полиэдра в виде сечения неотрицательного ортанта аффинным подпространством. Проблема формулируется в двух вариантах: в первом варианте требуется, чтобы размерность ортанта равнялась числу гиперграней полиэдра, а во втором ограничений на размерность ортанта не накладывается.

В первом варианте формулировки в настоящей работе проблема исследована до конца только в размерности два: мы приводим полную классификацию двумерных полиэдров, которые могут быть реализованы соответствующим образом. В старших размерностях по этой проблеме в работе получены частичные результаты. Что касается второго варианта формулировки проблемы, то в работе даётся почти исчерпывающий критерий о возможности реализации полиэдров таким способом (пробел касается лишь узкого класса неограниченных полиэдров). В частности, мы доказываем, что любой многогранник может быть реализован как сечение неотрицательного ортанта достаточно большой размерности.
\end{abstract}
%
%
\maketitle
%
%
\section{Введение}

Хорошо известно, что любой выпуклый $n$-мерный многогранник $P$ с точностью до аффинной эквивалентности может быть реализован как пересечение $$P=\alpha\cap\mathbbm{R}^m_{\geq 0}\subset \mathbbm{R}^m$$ $k$-мерной плоскости $\alpha$ и неотрицательного ортанта $\mathbbm{R}^m_{\geq 0}$, где $m$ равно числу гиперграней $P$ (см. \cite[Construction~1.2.1]{BP}, а также предложение~\ref{paff} настоящей работы). В то же время известный результат из школьного курса стереометрии гласит, что треугольник может быть получен как плоское сечение неотрицательного октанта $\mathbbm{R}^3_{\geq 0}$ в том и только том случае, если является остроугольным \cite{mccme}. Разумеется, тут речь идёт уже не об аффинной эквивалентности, а об изометрии. Эти два факта являются мотивировкой к постановке вопроса о том, какие многогранники или, более общо, какие полиэдры (о разнице в этих терминах см. раздел~\ref{intro}) можно реализовать в виде пересечения плоскости с неотрицательным ортантом. Причём вопрос можно поставить в двух вариантах.

\begin{question}\label{mainq}
Какие выпуклые $n$-мерные полиэдры с точностью до изометрии могут быть реализованы как пересечение $n$-мерной плоскости с положительным ортантом $\mathbbm{R}^m_{\geq 0}$, где $m$
\begin{itemize}
\item[А.] равно числу гиперграней полиэдра;
\item[B.] произвольно?
\end{itemize}
\end{question}

Вопрос A оказывается более интересным и сложным. Полиэдры, для которых ответ на вопрос А положительный, мы называем ортантными. Изучению свойства ортантности, главным образом, и посвящена настоящая работа. Наши основные результаты таковы:
\begin{itemize}
\item[1.] В теореме~\ref{t2d} мы приводим полную классификацию ортантных полиэдров размерности 2.
\item[2.] В теореме~\ref{too} мы сводим классификацию ортантных полиэдров произвольной размерности к классификации ортантных полиэдров с небольшим (в некотором смысле минимальным) числом граней, которые мы называем осн$\acute{\text{o}}$вными ортантными (см. определение~\ref{doo}).
\item[3.] В теореме~\ref{ts} мы приводим классификацию ортантных симплексов, а именно показываем, что ортантными среди симплексов являются в точности остроугольные ортоцентрические. Вместе с результатами теоремы~\ref{too} эта классификация позволяет определять ортантность более широкого класса полиэдров, а именно полиэдров, у которых ранг (см. определение~\ref{dr}) на единицу больше размерности.
\end{itemize}

Что же касается вопроса B, то в теореме~\ref{tqb} мы доказываем, что ответ на него очень простой: таким способом могут быть реализованы все многогранники (и все полиэдры с естественным ограничением на конус рецессии). Чтобы это доказать, однако, мы активно используем технику, развитую в ходе исследования вопроса A.

\medskip

Работа построена следующим образом.

В разделе \ref{intro} мы фиксируем обозначения и напоминаем о том, как полиэдры реализуются сечением ортанта с точностью до аффинной эквивалентности. В конце раздела мы кратко напоминаем имеющуюся теорию о вполне положительных матрицах, которая, как будет видно, напрямую связана с вопросом B для неограниченных полиэдров. Эта связь, однако, не будет изучаться в настоящей работе и приведена лишь для полноты изложения.

Разделы~\ref{rsist!}---\ref{rs} посвящены изучению вопроса A:
\begin{itemize}
 \item В разделе~\ref{rsist!} мы сводим вопрос ортантности полиэдра к вопросу существования положительного решения некоторой системы (!), строящейся по полиэдру.

\item В разделе~\ref{rez} мы делаем некоторые наблюдения о виде системы (!), которые позволяют  ввести на множестве полиэдров такое отношение эквивалентности, что вопросы ортантности эквивалентных полиэдров легко сводятся друг к другу. Факторизуя множество полиэдров по этому отношению эквивалентности, мы приходим к понятию ёжика полиэдра, с которым далее и работаем.

\item В разделе~\ref{r2d} мы классифицируем двумерные ортантные полиэдры методом грубой силы. При желании этот раздел может быть при чтении опущен, а классификация двумерных ортантных полиэдров выведена из результатов следующего раздела. Однако нам кажется, что этот раздел полезен для лучшего понимания происходящего.

\item В разделе~\ref{rgc} приводятся общие наблюдения, касающиеся свойства ортантности полиэдров произвольной размерности. Мы определяем ранг полиэдра, как ранг матрицы линейной системы (!), и показываем, что вопрос классификации ортантных полиэдров сводится к классификации ортантных полиэдров с минимальным числом граней при заданном ранге.

\item В разделе~\ref{rs} мы приводим классификацию ортантных полиэдров, у которых ранг на единицу больше размерности. С помощью результатов предыдущего раздела эта классификация сводится к классификации ортантных симплексов.
\end{itemize}

Раздел~\ref{rqb} посвящён ответу на вопрос~B.

Наконец, в разделе~\ref{rq} мы приводим вопросы, которые в настоящей работе остались не рассмотренными. Ответы на них автору не известны.

\section{Обозначения и предварительные наблюдения}\label{intro}

Основными объектами в работе являются выпуклые полиэдры. Разница в терминах ''полиэдр'' и ''многогранник'' для нас состоит в том, что полиэдр в отличие от многогранника не обязан быть ограниченным.

\begin{definition}
\emph{Выпуклым} $n$\emph{-мерным полиэдром} называется фигура $P$ в евклидовом пространстве $\R^n$, заданная конечной системой линейных неравенств
$$
\left\{ \begin{array}{ccc}
a_{11}x_1+\ldots+a_{1n}x_n&\geq&b_1\\
\ldots \\
a_{m1}x_1+\ldots+a_{mn}x_n&\geq&b_m
\end{array}
\right. \eqno (\ast)
$$
 и не содержащаяся ни в какой гиперплоскости. Ограниченный выпуклый полиэдр называется \emph{выпуклым многогранником}.
\end{definition}

Говоря о полиэдрах, далее мы всегда будем опускать слово ''выпуклый''. Систему ($\ast$), соответствующую полиэдру $P$, мы всегда будем полагать \emph{минимальной} в том смысле, что из неё нельзя выкинуть ни одно из уравнений так, чтобы она по-прежнему задавала полиэдр $P$. В этом случае неравенства системы ($\ast$) взаимно однозначно соответствуют граням коразмерности 1 полиэдра $P$:
грань $P_i$, соответствующая неравенству $a_{i1}x_1+\ldots+a_{in}x_n\geq b_i$, задаётся в $P$ уравнением $a_{i1}x_1+\ldots+a_{in}x_n=b_i$.

В отношении граней полиэдров мы будем придерживаться следующей терминологии: грани размерностей 0 и 1 мы называем, как обычно, вершинами и рёбрами, а грани коразмерности 1 мы называем гипергранями. Множество полиэдров размерности $n$ с $m$ гипергранями мы обозначаем $\Pol_n(m)$ (и опускаем число в скобках и (или) индекс в случае, если количество гиперграней и (или), соответственно, размерность не фиксированы). Аффинные подпространства размерности $k$ в $\R^n$ мы называем гиперплоскостями, если $k=n-1$, прямыми, если $k=1$, и плоскостями в случае произвольного $k$.

Иногда систему ($\ast$) мы будем записывать в матричном виде $A\textbf{x}\geq\textbf{b}$, а строки матрицы $A$ будем обозначать $\textbf{a}_i=(a_{i1},\ldots,a_{in})$. Полиэдр $P$ будем называть \emph{вырожденным}, если он содержит прямую. Легко заметить, что полиэдр $P$ размерности $n$ является невырожденным тогда и только тогда, когда соответствующая ему матрица $A$ имеет ранг $n$.

\begin{proposition}\label{paff}Пусть $P\in \Pol_n(m)$  --- невырожденный полиэдр, заданный системой ($\ast$) в $\R^n$. Тогда найдётся такая $n$-мерная плоскость $\alpha\subset \R^m$, что полиэдр $Q:=\alpha\cap \R^m_{\geq 0}$ является аффинно эквивалентным полиэдру $P$.
\end{proposition}
\begin{proof} Линейное отображение $\phi: \R^n \to \R^m$ $$(x_1,\ldots,x_n)\mapsto (a_{11}x_1+\ldots+a_{1n}x_n-b_1,\ldots,a_{m1}x_1+\ldots+a_{mn}x_n-b_m)$$ задаёт аффинную эквивалентность полиэдров $P$ и $\Im\phi\cap \R^m_{\geq 0}$.
\end{proof}
\begin{remark}Очевидно, что вырожденный полиэдр не может содержаться в неотрицательном ортанте, так что дальше в работе мы будем рассматривать только невырожденные полиэдры, опуская при этом слово ''невырожденный''.\end{remark}
\begin{remark}\label{rvideqv}Легко видеть, что отображение $\phi: \R^n \to \R^k$ задаёт аффинную эквивалентность полиэдров $P\subset \R^n$ и $Q=\Im \phi \cap \R^k_{\geq 0}$ тогда и только тогда, когда имеет вид
$$\phi(x_1,\ldots,x_n)=(u_{11}x_1+\ldots+u_{1n}x_n-v_1,\ldots,u_{k1}x_1+\ldots+u_{kn}x_n-v_k),$$
где все функции $f_i(\textbf{x}):=u_{i1}x_1+\ldots+u_{in}x_n-v_i$ неотрицательны на точках полиэдра $P$, причём для каждой гиперграни $P_i$ полиэдра $P$ найдётся ненулевая функция $f_j$, которая на этой грани зануляется, а значит имеет вид $f_j(\textbf{x})=k_i(a_{i1}x_1+\ldots+a_{in}x_n-b_{i})$ для некоторого $k_i\in \R_{>0}$. В частности, в случае $k=m$ такое отображение $\phi$ непременно имеет вид
$$\phi(x_1,\ldots,x_n)=$$  $$=\left(k_{\sigma(1)}(a_{\sigma(1)1}x_1+\ldots+a_{\sigma(1)n}x_n-b_{\sigma(1)}),\ldots,\right.$$ $$\left.k_{\sigma(m)}(a_{\sigma(m)1}x_1+\ldots+a_{\sigma(m)n}x_n-b_{\sigma(m)})\right),$$
где $\sigma\in \S_m$ --- некоторая перестановка.
\end{remark}

В дальнейшем полиэдры нас будут интересовать с точностью до изометрии.
\begin{definition} Полиэдры $P$ и $Q$ будем называть \emph{изометричными} и обозначать $P\cong Q$, если существует взаимно однозначное отображение $\phi: P\to Q$, сохраняющее расстояния между точками, то есть $\rho(\phi(\textbf{x}),\phi(\textbf{y}))=\rho(\textbf{x},\textbf{y})$ для любых $\textbf{x},\textbf{y}\in P$.\end{definition}

Вернёмся теперь к вопросам из введения, которым посвящена эта работа.
\begin{definition}Полиэдр $P\in\Pol_n(m)$ будем называть \emph{ортантным}, если существует такая ($n$-мерная) плоскость $\alpha\subset \R^m$, что $P\cong \alpha\cap \mathbbm{R}^m_{\geq 0}$, то есть $P$ может быть реализован так, как указано в вопресе~\ref{mainq}.A.\end{definition}
Что же касается полиэдров, удовлетворяющих вопросу~\ref{mainq}.B, то для них мы не будем вводить никакого специального термина по той простой причине, что, как уже говорилось, таковыми являются все полиэдры с одним естественным ограничением, касающимся только неограниченных полиэдров. Чтобы указать это ограничение, напомним следующее определение.
\begin{definition}\emph{Конусом рецессии} $C(P)$ полиэдра $P\subseteq \mathbbm{R}^n$ называется множество таких векторов $v\in\mathbbm{R}^n$, что  $x+v\in P$ для любой точки $x\in P$.
\end{definition}
\begin{remark}Необходимо отметить, что в этой работе конуса мы всегда рассматриваем как подмножества не аффинного евклидова пространства, а ассоциированного векторного.
\end{remark}
Если полиэдр $P$ оказался реализован как пересечение плоскости с неотрицательным ортантом $\mathbbm{R}^k_{\geq 0}$, то при этой реализации его конус рецессии $C(P)$ автоматически оказывается вложенным в $\mathbbm{R}^k_{\geq 0}$. В разделе~\ref{rqb} мы покажем, что при небольших дополнительных ограничениях верно и обратное: полиэдр $P$ может быть реализован как пересечение плоскости с неотрицательным ортантом, если $C(P)\setminus \{ \textbf{0} \}$  можно вложить в положительный ортант $\mathbbm{R}^k_{> 0}$ для некоторого $k$. В случае, когда $C(P)$ в неотрицательный ортант вкладывается, а в положительный --- нет, возникают неинтересные технические сложности, поэтому мы его опускаем.

\medskip

Определить по заданному полиэдральному конусу, можно ли его вложить в неотрицательный ортант, зачастую бывает сложно. Никакого общего алгоритма на этот счёт не известно. Маломерная интуиция тут легко может подвести. Так, возможность вложения конуса в неотрицательный ортант той же размерности не является необходимым условием. Простейший пример доставляет конус $\R_{\geq 0}\langle(1, 1, \sqrt{2}),(1,-1,\sqrt{2}),(-1,-1,\sqrt{2}),(-1,1,\sqrt{2})\rangle$, который не может быть вложен в $\R^3_{\geq 0}$, но изометричен конусу $\R^4_{\geq 0}\cap\{x_1+x_2-x_3-x_4=0\}$. Очевидным необходимым условием является то, что угол между любыми двумя векторами конуса не должен превосходить $\frac{\pi}{2}$. Более содержательное наблюдение состоит в том, что последнее условие не является достаточным!

Для полноты изложения вкратце напомним имеющуюся на этот счёт теорию, а за подробностями мы отсылаем в \cite{H}.

 Пусть $\textbf{v}_1,\ldots,\textbf{v}_m\in\R^n$ --- набор векторов. Рассмотрим их матрицу Грама
 $$\Gram[\textbf{v}_1,\ldots,\textbf{v}_m]:= \begin{pmatrix}
(\textbf{v}_1,\textbf{v}_1)& \ldots & (\textbf{v}_1,\textbf{v}_m)\\
\vdots& \ddots & \vdots\\
(\textbf{v}_m,\textbf{v}_1)& \ldots & (\textbf{v}_m,\textbf{v}_m)
\end{pmatrix}.$$
Симметрические матрицы, которые получаются таким образом, называются \emph{положительно полуопределёнными}. В том случае, если все элементы положительно полуопределённой матрицы неотрицательны, она называется \emph{двояконеотрицательной} (doubly nonnegative). Симметрическая матрица $A\in \R^{m\times m}$ называется \emph{вполне положительной} (completely positive), если она может быть представлена в виде $A=BB^T$, где $B\in \R^{m\times k}_{\geq 0}$.
\begin{proposition}Множества $\mathcal{DNN}_m$ двояконеотрицательных матриц и $\mathcal{CP}_m$ вполне положительных матриц образуют замкнутые выпуклые конусы в пространстве симметричных матриц размера $m\times m$, причём $\mathcal{CP}_m\subset \mathcal{DNN}_m$. Это вложение является строгим при $m\geq5$.
\end{proposition}

Теперь, если $\textbf{v}_1,\ldots,\textbf{v}_m$ --- это набор векторов на рёбрах полиэдрального конуса $C=\R_{\geq 0}\langle\textbf{v}_1,\ldots,\textbf{v}_m\rangle$, то, как несложно заметить, $\Gram[\textbf{v}_1,\ldots,\textbf{v}_m]$ является:
\begin{itemize}
\item двояконеотрицательной, если и только если угол между любыми двумя векторами из $C$ не превосходит $\frac{\pi}{2}$;
\item вполне положительной, если и только если $C$ может быть вложен в неотрицательный ортант.
\end{itemize}
Таким образом, понять, вкладывается ли заданный полиэдральный конус в неотрицательный ортант, --- это то же самое, что определить, является ли соответствующая симметрическая матрица вполне положительной. При этом необходимым, но не достаточным условием является то, что эта матрица является двояконеотрицательной. Общего алгоритма, позволяющего по двояконеотрицательной матрице определить, является ли она вполне положительной, никто не знает.

\section{Сведение вопроса ортантности к вопросу существования положительного решения некоторой линейной системы}\label{rsist!}

Пусть $P\in\Pol_n(m)$ --- полиэдр, заданный системой $(\ast)$ в $\R^n$. Мы хотим выяснить, можно ли его реализовать как пересечение $n$-мерной плоскости с неотрицательным ортантом в $\R^m$. Предположим, что это можно сделать. Тогда согласно замечанию \ref{rvideqv} соответствующая изометрия $\phi$ имеет следующий вид (по сравнению с тем общим видом, что указан в замечании~\ref{rvideqv}, мы опустили перестановку $\sigma$, так как она не влияет на свойство отображения $\phi$ являться изометрией):
$$(x_1,\ldots,x_n)\mapsto(k_1(a_{11}x_1+\ldots+a_{1n}x_n-b_{1}),\ldots,k_{m}(a_{m1}x_1+\ldots+a_{mn}x_n-b_{m})).$$
Вычислим $\rho(\phi(\textbf{x}),\phi(\textbf{y}))$, где $\textbf{x}=(x_1,\ldots,x_n)$, $\textbf{y}=(y_1,\ldots,y_n)$:
$$\rho\left(\phi(\textbf{x}),\phi(\textbf{y})\right)=\sqrt{\sum_{i=1}^{m} k_i^2\left(\sum_{j=1}^n a_{ij}(x_j-y_j)\right)^2}=$$
$$=\sqrt{\sum_{j=1}^n\sum_{i=1}^m k_i^2a_{ij}^2(x_j-y_j)^2+2\sum_{1\leq p<q\leq n}\sum_{r=1}^m k_r^2a_{rp}a_{rq}(x_p-y_p)(x_q-y_q)}.$$
Выполнение равенства $\rho(\phi(\textbf{x}),\phi(\textbf{y}))=\rho(\textbf{x},\textbf{y})$ для всех $\textbf{x},\textbf{y}\in P$ эквивалентно тому, что система из $C_{n+1}^2$ уравнений

$$
\left\{ \begin{array}{ccc}
a_{11}^2t_1+\ldots+a_{m1}^2t_m&=&1\\
\ldots \\
a_{1j}^2t_1+\ldots+a_{mj}^2t_m&=&1\\
\ldots \\
a_{1n}^2t_1+\ldots+a_{mn}^2t_m&=&1\\
a_{11}a_{12}t_1+\ldots+a_{m1}a_{m2}t_m&=&0\\
\ldots\\
a_{1p}a_{1q}t_1+\ldots+a_{mp}a_{mq}t_m&=&0\\
\ldots\\
a_{1,n-1}a_{1n}t_1+\ldots+a_{m,n-1}a_{mn}t_m&=&0
\end{array}
\right. \eqno (!)
$$
имеет положительное решение $\textbf{t}=(t_1,\ldots,t_m)=(k_1^2,\ldots,k_m^2)\in\R^m_{>0}$.

Мы доказали следующее утверждение.
\begin{theorem}\label{tsist} Полиэдр $P\in\Pol_n(m)$, заданный системой неравенств $(\ast)$, является ортантным тогда и только тогда, когда система уравнений $(!)$ имеет положительное решение.
\end{theorem}

\section{Ёжик полиэдра}\label{rez}

Идея этого раздела предельно проста. В предыдущем разделе мы показали, что свойство ортантности полиэдра зависит только от свойства соответствующей ему системы $(!)$ иметь положительное решение. Теперь, если вопрос существования положительного решения системы, соответствующей одному полиэдру, легко сводится к вопросу существования положительного решения системы, соответствующей другому полиэдру, то нет и никакого смысла исследовать проблему ортантности для этих полиэдров по отдельности. Исходя из этого, мы хотим выделить как можно более большие классы полиэдров, для которых вопросы существования положительного решения соответствующих систем $(!)$ сводятся один к другому, и назвать такие полиэдры эквивалентными. После факторизации по этому отношению эквивалентности от полиэдра остаются данные, по которым система $(!)$ восстанавливается с точностью до простых преобразований, не меняющих её свойства иметь положительное решение. Мы называем эти данные ёжиком. Ёжики определяются меньшим числом параметров, чем система $(\ast)$, задающая полиэдр, по которой в предыдущем разделе мы строили систему $(!)$. Свойство ортантности полиэдров спускается на ёжики и, учитывая вышесказанное, именно для них проблему ортантности разумно будет далее изучать.

\medskip

Пусть $P\in\Pol_n(m)$ задан системой $(\ast)$, по которой затем, как в разделе~\ref{rsist!}, построена система $(!)$. Выделим те преобразования, которые можно производить с системой $(\ast)$, чтобы свойство системы $(!)$ иметь положительное решение не менялось.
\begin{itemize}
\item[I.] Можно как угодно изменять столбец свободных членов $\textbf{b}$: от него вид системы $(!)$ вообще не зависит.
\item[II.] Можно умножать векторы $\textbf{a}_i$ на произвольные ненулевые скаляры. Такое преобразование соответствует умножению соответствующего столбца матрицы $(!)$ на положительное число, что не меняет свойства иметь положительное решение.
\item[III.] Можно умножать $A$ на произвольную ортогональную матрицу из $\O(n)$. После этого новая система $(!)$ будет иметь такие же решения, как и старая. И не удивительно: по модулю уже отмеченного преобразования I умножение $A$ на ортогональную матрицу задаёт поворот полиэдра, который не может изменить его свойства ортантности.
\item[IV.] Можно добавлять новые неравенства вида $\textbf{a}_{new}\textbf{x}\geq b_{new}$, где $\textbf{a}_{new}=k\textbf{a}_i$ для некоторых $i$ и $k\ne 0$, а $b_{new}$ --- произвольное.
\item[V.] Если $\textbf{a}_j=k\textbf{a}_i$ для некоторого $k\ne 0$ и при этом $i\ne j$, то одно из двух неравенств, соответствующих этим векторам, можно выбросить.
\end{itemize}

\begin{definition} Полиэдры $P$ и $Q$ будем называть \emph{эквивалентными} и обозначать $P\sim Q$, если системы $(\ast)$, которые их задают, могут быть получены одна из другой последовательностью преобразований вида I-V. \emph{Ёжиками} будем называть элементы множества $\Y:=\Pol/\sim$, то есть классы эквивалентности полиэдров. Обозначим $\yz: \Pol\to \Y$ отображение факторизации. Элемент $\yz(P)\in\Y$, будем называть \emph{ёжиком} полиэдра $P$.\end{definition}

Отметим, что имеется разбиение $\Y=\coprod \Y_{n}(m)$, где $\y\in\Y_{n}(m)$, если $m$ --- это минимум количества гиперграней полиэдров, находящихся в этом классе эквивалентности (в дальнейшем мы будем говорить, что $m$ --- это \emph{число иголок} ёжика $\y$), а $n$ --- это размерность полиэдров в этом классе (она, естественно, у всех представителей одинакова).


\medskip

Пусть $P\in\Pol_n(m)$, и пусть максимальное число непараллельных граней, которые можно выбрать в $P$, равно $m'$. Отметив на единичной сфере концы векторов нормали произвольного набора из $m'$ непараллельных гиперграней, полиэдру $P$ можно корректно сопоставить элемент $\widetilde{\yz}(P)$ множества $$\widetilde{\Y}_{n}(m'):=\left(\left.\left(\left.S^{n-1}\right/(\Z/2\Z)\right)^{m'}\setminus\Delta\right)\right/\left(\S_{m'}\times \O(n)\right),$$ где $S^{n-1}\subset\R^n$ --- единичная сфера, на которой $\Z/2\Z$ действует центральной симметрией, $\Delta$ --- большая диагональ (те элементы произведения, которые совпадают при проекции на два разных сомножителя), $\S_{m'}$ действует на сомножителях перестановками, а $\O(n)$ действует симметриями одновременно на всех сомножителях. Далее можно заметить, что $$\widetilde{\yz}(P)=\widetilde{\yz}(Q)\text{ }\Longleftrightarrow \text{ }\yz(P)=\yz(Q),$$
поэтому $\Y_n(m')\subset \widetilde{\Y}_{n}(m')$.

Пространство $\left.S^{n-1}\right/(\Z/2\Z)$ можно представлять как верхнюю полусферу $S^{n-1}\cap\{x_n\geq 0\}$, у которой отождествлены противоположные точки границы $S^{n-2}\cong S^{n-1}\cap\{x_n= 0\}$.
   Ёжик $\y\in\Y_n(m')$ мы будем представлять и изображать как $m'$ векторов, проведённых из начала координат к точкам верхней полусферы (помня о том, что ёжики совпадают не только в том случае, когда совпадают их изображения). Векторы, из которых состоит изображение ёжика, будем называть его \emph{иголками}.

\begin{remark}На самом деле, $\widetilde{\Y}_{n}(m)=\Y_{n}(m)\sqcup\widetilde{\Y}_{n-1}(m)$. В самом деле, пусть $\widetilde{\y}\in \widetilde{\Y}_{n}(m)$. Поднимем $\widetilde{\y}$ произвольным образом до элемента в $(S^{n-1})^m$. Этот элемент соответствует набору из $m$ различных точек на сфере $S^{n-1}\subset\R^n$. Теперь, если эти точки лежат в одной гиперплоскости, то $\widetilde{\y}\in \widetilde{\Y}_{n-1}(m)$. Если же нет, то проведём касательные гиперплоскости к единичной сфере в этих точках. Эти гиперплоскости делят пространство на несколько полиэдральных областей. Обозначим $P$ ту область, что содержит центр сферы. Легко проверить, что $\widetilde{\y}=\widetilde{\yz}(P)$.
\end{remark}

Зафиксируем, что ортантность полиэдра корректно спускается на множество ёжиков.
\begin{proposition}\label{yozh}Пусть $P\sim Q$. Тогда полиэдр $P$ является ортантным в том и только том случае, если ортантным является полиэдр $Q$.\end{proposition}
\begin{proof}Сразу следует из определения эквивалентных полиэдров и теоремы~\ref{tsist}.
\end{proof}

\begin{definition}Ёжик $\y\in\Y$ называется \emph{ортантным}, если $\y=\yz(P)$, где $P$ --- ортантный полиэдр.\end{definition}

На множество ёжиков спускаются и другие определения, связанные со свойством ортантности, например, ранг полиэдра и свойство полиэдра быть осн$\acute{\text{o}}$вным ортантным, которые мы введём в разделе~\ref{rgc}.

\begin{example}Трапеция, изображённая на рисунке ниже, не является ортантной, поскольку её ёжик совпадает с ёжиком прямоугольного треугольника, про который мы знаем, что он не ортантный. А вот правильный шестиугольник является ортантным, так как его ёжик совпадает с ёжиком равностороннего треугольника.

\newdimen\Rad
\Rad=0.8cm
\begin{center}
\begin{tikzpicture}
    \draw[yshift=-0.35\Rad] (0,0) -- (0,1.2\Rad) -- (1.8\Rad,-0.6\Rad) -- (0.6\Rad,-0.6\Rad) -- (0,0) node at (1.9\Rad,0.35\Rad) {$\sim$};
    \draw[xshift=2.7\Rad, yshift=-0.8\Rad] (0,0) -- (0,1.6\Rad) -- (1.6\Rad,0) -- (0,0);

    \draw[xshift=8.5\Rad] (30:\Rad) \foreach \x in {30,90,...,330} {
            -- (\x:\Rad)
    }-- cycle (0:1.5\Rad) node {$\sim$};
    \draw[xshift=11.5\Rad] (60:\Rad) \foreach \x in {180,300} {
            -- (\x:\Rad)
    } -- cycle;
    \draw[xshift=2\Rad,yshift=-2.3\Rad,thick,lightgray] ([shift=(0:\Rad)]0:0) arc (0:180:\Rad);
    \draw[xshift=9.85\Rad,yshift=-2.3\Rad,thick,lightgray] ([shift=(0:\Rad)]0:0) arc (0:180:\Rad);
    \foreach \x in {0,45,90} {
    \draw[xshift=2\Rad,yshift=-2.3\Rad, thick, ->]
       (0:0)--(\x:\Rad);
    }
     \foreach \x in {0,60,120} {
    \draw[xshift=9.85\Rad,yshift=-2.3\Rad, thick, ->]
       (0:0)--(\x:\Rad);
    }

\end{tikzpicture}
\end{center}

\end{example}

\begin{definition}Две системы линейных уравнений будем называть \emph{эквивалентными}, если либо обе они имеют положительное решение, либо ни одна из них положительного решения не имеет.\end{definition}
\begin{corollary}Если $P\sim Q$, то соответствующие им системы $(!)$ эквивалентны.\end{corollary}

\begin{definition}
Будем говорить, что ёжик $\y'$ является \emph{подъёжиком} в ёжике $\y$, если их можно изобразить одновременно так, что иголки $\y'$ будут являться подмножеством иголок $\y$. В этом случае будем писать $\y' \subset \y$ и $P\prec Q$, если $\y'=\yz(P)$, $\y=\yz(Q)$.
\end{definition}

\begin{remark}\label{ziez}Ясно, что $\prec$ является отношением частичного порядка. В явном виде $P\prec Q$ означает, что эти полиэдры могут быть так повёрнуты, что для каждой гиперграни полиэдра $P$ найдётся параллельная ей гипергрань $Q$. При этом $P\simeq Q$ тогда и только тогда, когда $P\prec Q$ и $Q\prec P$.\end{remark}

\begin{remark}Пусть полиэдр $P$ задан системой $(\ast)$. Выкинем несколько неравенств из этой системы и рассмотрим полиэдр $Q$, который задают оставшиеся неравенства. Этот полиэдр образуется, как область пространства, ограниченная гиперплоскостями, на которых лежат гиперграни $P$, соответствующие оставшимся неравенствам системы. Легко заметить, что $Q\prec P$ и, наоборот, все $\y'\subset \yz(P)$ реализуются таким образом.
\end{remark}

Вообще говоря, свойство $P\prec Q$ не позволяет судить об ортантности $P$ по ортантности $Q$ и наоборот, о чём свидетельствует пример ниже. Тем не менее, следствие~\ref{coo} позволит нам сделать некоторые важные наблюдения по этому поводу.
\begin{example}
Как легко заметить,

\begin{center}
\begin{tikzpicture}
\newdimen\J
\J=1cm
    \draw (45:\J) \foreach \x in {45,135,...,315} {
            -- (\x:\J)
    }-- cycle (0:1.5\J) node {$\prec$};

    \draw[xshift=2.2\J, yshift=-0.8\J] (0,0) -- (0,1.6\J) -- (1.6\J,0) -- (0,0) node at (2\J,0.8\J) {$\prec$};

    \draw[xshift=6\J] (22.5:\J) \foreach \x in {22.5,67.5,...,337.5} {
            -- (\x:\J)
    } -- cycle (0:\J) node[below right]{$,$};
\end{tikzpicture}
\end{center}
однако, как будет ясно из результатов раздела~\ref{r2d}, квадрат и правильный восьмиугольник являются ортантными, а прямоугольный треугольник --- нет.
\end{example}

Далее, чтобы работать с ёжиками, нам часто будет необходимо возвращаться от ёжика $\y$ обратно к полиэдру и сопоставлять ёжику системы $(\ast)$ и $(!)$. Никакого канонического способа выбрать полиэдр в $\yz^{-1}(\y)$ нет, но, тем не менее, в $\yz^{-1}(\y)$ есть полиэдры, которые вместе с соответствующими системами $(\ast)$ выглядят несколько проще остальных.

Пусть сначала $P\in\yz^{-1}(\y)$ --- произвольный, $(\ast)$ --- соответствующая ему система. Обозначим $\textbf{a}_j:=(a_{j1},\ldots,a_{jn})$. Зафиксируем среди этих векторов упорядоченный набор из $n-1$ линейно независимых. Можно считать, что это векторы $\textbf{a}_1,\ldots,\textbf{a}_{n-1}$. По сделанному выбору можно каноническим образом построить полиэдр $Q$ эквивалентный $P$ с более простым видом систем $(\ast)$ и $(!)$.

Проделаем это в несколько шагов.
\begin{itemize}
\item Произведя нужный поворот полиэдра $P$, можно считать, что векторы $\textbf{a}_j$ для $1\leq j\leq n-1$ имеют вид $\textbf{a}_1=(a_{11},0,\ldots,0),\,\textbf{a}_2=(a_{21},a_{22},0,\ldots,0),\,\textbf{a}_3=(a_{31},a_{32},a_{33},0,\ldots,0),\,\ldots,\,\textbf{a}_{n-1}=(a_{n-1,1},a_{n-1,2},\ldots,a_{n-1,n-1},0)$.
\item Поделим все векторы $\textbf{a}_j$ на их норму и теперь будем полагать $\|\textbf{a}_j\|=1$. После соответствующей корректировки чисел $b_j$, система $(\ast)$ будет задавать всё тот же полиэдр.
\item Пусть теперь $1\leq j\leq m$ и $\textbf{a}_j=(a_{j1},a_{j2},\ldots,a_{jp_j},0,\ldots,0)$, где $a_{jp_j}\ne 0$. Тогда, умножив при необходимости вектор $\textbf{a}_j$ на $-1$, можно считать, что $a_{jp_j}>0$. Если в полученном наборе векторов оказались равные, то оставим из них ровно по одному представителю, сохранив векторы $\textbf{a}_1$, ..., $\textbf{a}_{n-1}$. Положим $b_j=-1$ для всех $j$.
\end{itemize}

Легко видеть, что после последнего шага система $(\ast)$ минимальна и задаёт полиэдр $Q\sim P$, а потому $Q\in\yz^{-1}(\y)$. Эту систему мы будем называть \emph{приведённым видом} исходной системы, а $Q$ --- \emph{приведённым полиэдром}. Изображение полученных векторов $\textbf{a}_j$, отложенных от начала координат, будем называть \emph{приведённым изображением} ёжика $\y$. Сопоставляя каждому ёжику произвольный приведённый полиэдр $s(\y)\in \yz^{-1}(\y)$, зафиксируем сечение $s:\Y\to \Pol$. Отметим, что $s(\Y_n(m))\subset \Pol_n(m)$. Всегда будем подразумевать, что система $(\ast)$, соответствующая $s(\y)$, имеет приведённый вид.

\section{Классификация ортантных полиэдров размерности 2}\label{r2d}

Пусть $\y\in\Y_n(m)$, $s(\y)$ --- приведённый двумерный полиэдр. Для него система $(!)$ состоит всего из трёх уравнений. После перенумерации можно записать $\textbf{a}_j=(\cos\psi_j,\sin\psi_j)$, где $0=\psi_1<\psi_2<\ldots<\psi_m<\pi$.
Система $(!)$ в таком случае принимает вид:
$$
\left\{ \begin{array}{ccl}
(\cos^2\psi_1)t_1+\ldots+(\cos^2\psi_m) t_m&=&1\\
(\sin^2\psi_1)t_1+\ldots+(\sin^2\psi_m) t_m&=&1\\
(\cos\psi_1\sin\psi_1)t_1+\ldots+(\cos\psi_m\sin\psi_m)t_m&=&0.
\end{array}
\right.\eqno (1)
$$
Вычтя из первого уравнения системы второе, получим:
$$
\left\{ \begin{array}{ccl}
(\cos^2\psi_1-\sin^2\psi_1)t_1+\ldots+(\cos^2\psi_m-\sin^2\psi_m)t_m&=&0\\
(\cos\psi_1\sin\psi_1)t_1+\ldots+(\cos\psi_m\sin\psi_m)t_m&=&0.
\end{array}
\right.\eqno (2)
$$
Заметим, что система (2) эквивалентна системе (1): если $(t_1,\ldots,t_m)$ --- положительное решение системы (2), то $\left(\frac{t_1}{c},\ldots,\frac{t_m}{c}\right)$, где $c=\sum (\cos^2\psi_j) t_j$, --- положительное решение системы (1).

Воспользовавшись формулой двойного угла, получим:
$$
\left\{ \begin{array}{ccl}
(\cos2\psi_1)t_1+\ldots+(\cos2\psi_m)t_m&=&0\\
(\sin2\psi_1)t_1+\ldots+(\sin2\psi_m)t_m&=&0.
\end{array}
\right.\eqno (3)
$$
Система (3) имеет положительное решение, тогда и только тогда, когда конус $\R_{\geq 0} \langle (\cos2\psi_1,\sin2\psi_1),\ldots,(\cos2\psi_m,\sin2\psi_m)\rangle$ является линейным подпространством в $\R^2$ (делая подобные наблюдения, часто упоминают лемму Фаркаша). Отсюда следует, что положительное решение найдётся в точности в следующих случаях:
\begin{itemize}
\item[А:] $m=2$, $\psi_2=\frac{\pi}{2}$;
\item[Б:] $m>2$, $\psi_m>\frac{\pi}{2}$ и $\psi_{i+1}-\psi_i<\frac{\pi}{2}$ для любого $i\leq m-1$.
\end{itemize}
Вариант А мы называем \emph{вырожденным случаем}, а вариант Б --- \emph{невырожденным случаем}. Точным смыслом мы наделим эти понятия в разделе~\ref{rgc}. Как несложно заметить, условие Б невырожденного случая можно переписать в виде
$$m>2,\,\psi_m>\frac{\pi}{2}\text{ и }\psi_{p+1}-\psi_p<\frac{\pi}{2}\text{, где }p=\max_{\psi_i<\frac{\pi}{2}} i.\eqno (\maltese)$$

\begin{remark}\label{zot}Пусть $\y\in \Y_2(3)$. Рассмотрим произвольное изображение $\y$ в верхнем полукруге:

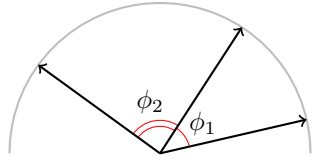
\begin{center}
\begin{tikzpicture}
\newdimen\J
\J=2cm
    \draw[thick,lightgray] ([shift=(0:\J)]0:0) arc (0:180:\J);
    \draw[red] ([shift=(13:0.2\J)]0:0) arc (13:57:0.2\J);
       \draw[red] ([shift=(57:0.22\J)]0:0) arc (57:144:0.22\J);
       \draw[red] ([shift=(57:0.18\J)]0:0) arc (57:144:0.18\J);
    \foreach \x in {13,57,144} {
    \draw[thick, ->]
       (0:0)--(\x:\J);
    }

       \draw node at (35:0.35\J) {$\phi_1$};
       \draw node at (100.5:0.35\J) {$\phi_2$};

\end{tikzpicture}
\end{center}
В $\yz^{-1}(\y)$ найдётся остроугольный треугольник тогда и только тогда, когда $\phi_1<\frac{\pi}{2}$, $\phi_2<\frac{\pi}{2}$, а $\phi_1+\phi_2>\frac{\pi}{2}$.
\end{remark}

Это замечание позволяет переформулировать условие для невырожденного случая в следующем виде.

\begin{proposition}\label{p2d}Пусть $\y\in\Y_n(m)$, $m>2$. Тогда $\y$ является ортантным в точности в двух следующих (непересекающихся) случаях:
\begin{itemize}
\item[(Б1)] $m=4$, $\y$ состоит из двух пар перпендикулярных иголок, то есть имеет следующее изображение:
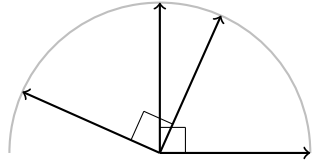
\begin{center}
\begin{tikzpicture}
\newdimen\J
\J=2cm
    \draw[thick,lightgray] ([shift=(0:\J)]0:0) arc (0:180:\J);
    \foreach \x in {0,66,90,156} {
    \draw[thick, ->]
       (0:0)--(\x:\J);
    }
       \draw[shift=(90:0.17\J)] (0:0)-- (0:0.17\J);
       \draw[shift=(0:0.17\J)] (0:0)-- (90:0.17\J);
       \draw[shift=(156:0.21\J)] (0:0)-- (66:0.21\J);
       \draw[ shift=(66:0.21\J)] (0:0)-- (156:0.21\J);


\end{tikzpicture}
\end{center}
\item[(Б2)] существует такой остроугольный треугольник $T$, что $\y(T)\subset \y$.
\end{itemize}
\end{proposition}
\begin{proof}
$\underline{\Rightarrow}$ Пусть $\y$ --- ортантный. Так как $m>2$, мы в невырожденном случае и иголки $\y$ удовлетворяют условию $(\maltese)$. Обозначим $\y(i,j,k)$ --- подъёжик в $\y$, образованный иголками с номерами $i,j,k$. Пусть подъёжика, являющегося ёжиком остроугольного треугольника, в $\y$ нет. Тогда, применяя замечание~\ref{zot}
\begin{itemize}
\item[1.] к $\y(1,p,p+1)$, находим, что $\psi_{p+1}=\frac{\pi}{2}$;
\item[2.] к $\y(1,p,p+2)$, находим, что $\psi_{p+2}-\psi_p\geq\frac{\pi}{2}$;
\item[3.] к $\y(p,p+1,p+2)$, находим, что $\psi_{p+2}-\psi_p\leq\frac{\pi}{2}$, откуда $\psi_{p+2}-\psi_p=\frac{\pi}{2}$.
\end{itemize}
Теперь по замечанию~\ref{zot} ёжиком остроугольного треугольника является:
\begin{itemize}
\item[4.] подъёжик $\y(2,p+1,p+2)$, если $2<p$;
\item[5.] подъёжик $\y(p,p+1,m)$, если $m>p+2$.
\end{itemize}
Поэтому в $\y$ всего четыре иголки и он имеет тот вид, что указан в пункте Б1.

$\underline{\Leftarrow}$ Иголки канонических изображений ёжиков остроугольных треугольников удовлетворяют условию Б согласно замечанию~\ref{zot}. Остаётся заметить, что добавление новых иголок не может нарушить выполнение этого условия. Кроме того, простая проверка показывает, что ёжики вида Б1 также удовлетворяют условию Б.
\end{proof}
\begin{remark}\label{z2d}Условие Б1 можно было бы ослабить до того, что в $\y$ найдётся подъёжик, имеющий соответствующий вид, но тогда условия Б1 и Б2 могли бы выполняться одновременно. Полученный результат, на самом деле, несложно вывести из общей теоремы~\ref{too} из следующего раздела работы.
\end{remark}

Результаты этого раздела скомпонуем следующим образом (тут эквивалентность $(i)\Leftrightarrow (iii)$ --- это просто переписанное с учётом замечаний~\ref{z2d} и~\ref{ziez} предложение~\ref{p2d}).
\begin{theorem}\label{t2d}
Следующие условия эквивалентны:
\begin{itemize}
\item[(i)] Двумерный полиэдр $P$ является ортантным;
\item[(ii)] Каноническое изображение ёжика полиэдра $P$ состоит из двух перпендикулярных иголок или удовлетворяет условию $(\maltese)$;
\item[(iii)] Выполнено хотя бы одно из трёх условий:
\begin{itemize}\item в полиэдре $P$ имеются такие три стороны, что прямые, на которых они лежат, образуют в пересечении остроугольный треугольник;
\item в полиэдре $P$ найдутся две пары перпендикулярных сторон, что ни одна сторона одной пары не совпадает и не параллельна сторонам другой пары;
\item $P$ является одним из следующих полиэдров:
        \begin{itemize}
          \item[$\llcorner$] --- уголком;
          \item[$\sqsubset$] --- полосой с перпендикулярным границе разрезом;
          \item[$\square$] --- прямоугольником.
        \end{itemize}
     \end{itemize}
\end{itemize}
\end{theorem}
Теорема~\ref{t2d} позволяет моментально сделать следующее наблюдение.
\begin{corollary}
Все правильные многоугольники являются ортантными.
\end{corollary}

\section{Общие наблюдения, касающиеся ортантных полиэдров старших размерностей}\label{rgc}

Из предыдущего раздела следует отметить, что классификация двумерных ортантных полиэдров, поделилась на две части: то, что мы назвали невырожденным и вырожденным случаями. Такое же деление будет происходить и в старших размерностях, но из-за быстрого увеличения числа уравнений в системе $(!)$ с ростом размерности полиэдра, число вырожденных случаев будет стремительно расти, и исследовать каждый конкретный случай будет гораздо сложнее. Поясним, что под этим имеется в виду.

Пусть $P$ --- полиэдр размерности $n$, заданной системой $(\ast)$. Обозначим $Q$ матрицу левой части соответствующей системы уравнений $(!)$, а $Q'$ --- расширенную матрицу этой системы.
\begin{remark}Число $\rk Q$ является аффинным инвариантом полиэдра. Число $\rk Q'$ не является аффинным инвариантом полиэдра, но сохраняется при преобразованиях подобия. Кроме того, числа $\rk Q$ и $\rk Q'$ корректно спускаются на множество ёжиков. Если полиэдр (или ёжик) является ортантным, то $\rk Q=\rk Q'$. Обратное неверно.\end{remark}
\begin{definition}\label{dr}Число $\rk Q$ будем называть \emph{рангом} полиэдра (ёжика).\end{definition}

Будем обозначать $\Pol^r_n(m)\subset \Pol_n(m)$ и $\Y^r_n(m)\subset \Y_n(m)$ подмножества полиэдров и ёжиков ранга $r$ (и опускать в этих обозначениях скобки, если число гиперграней не фиксировано).

Задача классификации ортантных полиэдров размерности $n$ разбивается на подзадачи классификации ортантных полиэдров размерности $n$ и ранга $r$. К невырожденному случаю мы относим полиэдры максимального ранга $C^2_{n+1}$. В этом случае условие $\rk Q=\rk Q'$ не накладывает никаких дополнительных ограничений, и поэтому ортантность является открытым условием в пространстве таких полиэдров (этим словам несложно придать формальный смысл, но мы не будем этого делать). Если же ранг $r$ ортантного полиэдра $P$ строго меньше $C^2_{n+1}$, то этот полиэдр обязан удовлетворять метрическим соотношениям, соответствующим занулению всех миноров порядка $r+1$ матрицы $Q'$. Это условие замкнутое и только уже для полиэдров, удовлетворяющих этому условию, ортантность будет условием открытым.

\begin{example}\label{endgo}Пример ортантного полиэдра из $\Pol^{C^2_{n+1}}_n$ доставляет полиэдр, заданный неравенствами $x_i-x_j\geq -1,\, x_i+x_j\geq -1,\, x_i\geq -\frac{2}{3}$: несложно проверить, что система (!), соответствующая этому полиэдру, имеет положительное решение, а её ранг равен $C^2_{n+1}$.
\end{example}

Наименьший ранг, который может быть у (невырожденного) полиэдра размерности $n$, равен $n$. С классификацией ортантных полиэдров, у которых размерность равна рангу всё просто: ортантными среди них являются в точности те полиэдры, что эквивалентны полиэдру $\R^n_{\geq 0}$. Случаю полиэдров ранга $n+1$ посвящён следующий раздел работы. Для таких полиэдров тоже можно дать полный ответ. В размерности 2 всё ровно этими двумя случаями и ограничивается. В старших же размерностях основную сложность представляет именно классификация ортантных полиэдров ранга $r$ большего $n+1$. В настоящей работе классифицировать такие полиэдры мы даже и не пытаемся. Ниже мы лишь покажем, что эта проблема сводится к классификации ортантных полиэдров ранга $r$ с числом гиперграней $r$ (то есть минимально возможным) или, что тоже самое, ортантных ёжиков ранга $r$, у которых $r$ иголок.

\begin{definition}\label{doo}
Ортантные полиэдры из $\Pol_n^r(r)$ и ортантные ёжики из $\Y_n^r(r)$ называются \emph{осн}$\acute{o}${вными ортантными}.
\end{definition}

\begin{example}В размерности 2 основными ортантными полиэдрами являются уголок и остроугольные треугольники.
\end{example}

Основное утверждение этого раздела мы сформулируем для ёжиков, а читатель при желании может, используя замечание~\ref{ziez}, легко восстановить эквивалентную формулировку для полиэдров.
\begin{theorem}\label{too}
Ёжик $\y\in \Y_n^r(m)$ является ортантным тогда и только тогда, когда существуют основные ортантные ёжики $\y_1,\ldots,\y_k\subset\y$ такие, что $\y_1\cup\ldots\cup\y_k$ --- ёжик ранга $r$.
\end{theorem}

Объединение подъёжиков здесь понимается как ёжик, образованный объединением их иголок на некотором фиксированном изображении ёжика $\y$.

\begin{proof}
$\underline{\Rightarrow}$ Ведём индукцию по $n$ и $m$. Если $m=r$, доказывать нечего. Пусть $m>r$, а утверждение доказано для всех ёжиков размерности меньшей $n$ и всех ёжиков размерности $n$, у которых не более $m-1$ иголок. Положим $P=s(\y)$ и рассмотрим соответствующую систему $(!)$. Так как $\y$ ортантный, она имеет некоторое положительное решение $\textbf{t}$. При этом все решения этой системы образуют $(m-r)$-мерную плоскость в $\R^m$. Выберем в этой плоскости произвольную прямую, проходящую через $\textbf{t}$. Несложно заметить, что такая прямая пересекает множество $\R^m_{\geq 0}\setminus \R^m_{>0}$ в двух точках; обозначим эти точки $\textbf{u}=(u_1,\ldots,u_m)$ и $\textbf{v}=(v_1,\ldots,v_m)$. Обозначим $I=\{i:u_i=0\}$, $J=\{i:v_i=0\}$. Отметим, что оба эти множества непустые, а вот $I\cap J=\emptyset$. Рассмотрим полиэдры $P_I$ и $P_J$, получающиеся из системы $(\ast)$, задающей $P$, отбрасыванием неравенств с номерами из множеств $I$ и $J$ соответственно. Оба эти полиэдра являются ортантными, так как $\textbf{u}$ и $\textbf{v}$ доставляют положительные решения соответствующих систем $(!)$, и $\y(P_I)\cup\y(P_J)=\y$. Остаётся лишь воспользоваться предположением индукции для: \begin{itemize}
\item $\y(P_I)$ или $\y(P_J)$, если ранг одного из этих ёжиков равен $r$;
\item одновременно для $\y(P_I)$ и $\y(P_J)$, если ранг обоих ёжиков меньше $r$.
\end{itemize}
$\underline{\Leftarrow}$ Пусть $P=s(\y)$. По замечанию~\ref{ziez} каждый из ёжиков $\y_i$ соответствует выбору набора гиперграней полиэдра $P$ или, что то же самое, выбору подмножества $I_i\subset \{t_1,...,t_m\}$ переменных в системе уравнений $(!)$, соответствующей $P$. Нам надо доказать, что система $(!)$ имеет положительное решение. Так как $\y_i$ --- ортантный, система $(!)$ имеет решение $t_j=s_{ij}$, где $s_{ij}>0$, при $t_j\in I_i$ и $t_j=0$ иначе. Следовательно, система $(!)$ имеет решение
$$t_j=\frac{1}{k}\sum_{\{i: t_j\in I_i\}} s_{ij}.$$
Заметим, что в этом решении $t_j>0$, если $t_j\in \cup I_i$, и $t_j=0$ иначе. Поскольку ранг матрицы, образованной столбцами $\cup I_i$, равен рангу всей матрицы, соответствующей системе $(!)$, можно выбрать такое $\varepsilon >0$, что, положив $t_j=\varepsilon$ при $t_j\not\in \cup I_i$, система $(!)$ будет иметь положительное решение относительно $t_j\in \cup I_i$.
\end{proof}

\begin{corollary}\label{coo}Пусть $\y,\y'\in\Y^r_n$, $\y'\subset y$. Тогда, если $\y'$ ортантный, то и $\y$ является таковым.\end{corollary}

В частности, если для некоторого полиэдра $P\in\Pol_n$ мы видим, что гиперплоскости, на которых лежит некоторый набор его граней образуют в пересечении ортантный полиэдр ранга $C^2_{n+1}$, то $P$ является ортантным (и имеет тот же ранг).

\medskip

Раздел~\ref{r2d} мы завершили применением полученных результатов к правильным многоугольникам. В старших размерностях результатов мы получили мало, но тем не менее с сериями правильных многогранников несложно разобраться.
\begin{proposition}
Правильный n-мерный симплекс, n-мерный куб и n-мерный октаэдр являются ортантными.
\end{proposition}

\begin{proof}
Правильные симплексы не только попадают под теорему~\ref{ts}, которую мы докажем в следующем разделе, но и просто высекаются в ортанте $\R^{n+1}_{\geq 0}$ плоскостью $\{x_1+\ldots+x_{n+1}=0\}$. Куб размерности $n$ эквивалентен ортанту $\R^n_{\geq 0}$, который реализуется тривиальным образом. С октаэдром придётся немного повозиться.

Рассмотрим систему $(!)$, соответствующую $n$-мерному окдаэдру, заданному системой из $2^n$ неравенств $\pm x_1\pm\ldots\pm x_n\geq -1$. Первые $n$ уравнений этой системы одинаковы: $t_1+\ldots+t_{2^n}=1$. Про них можно забыть, так как по положительному решению оставшейся системы легко построить положительное решение всей системы. Рассмотрим матрицу $Q=(q_{ij})$ оставшейся однородной системы. Её строки нумеруются парами $(a,b)$, где $1\leq a<b\leq n$, а столбцы --- гипергранями октаэдра. Выберем те столбцы, что соответствуют гиперграням $\pm x_1\pm\ldots\pm x_n= -1$, где минусов ровно $\left[\frac{n}{2}\right]$. Пусть это последние $r:=C^{\left[\frac{n}{2}\right]}_n$ столбцов, а первый столбец соответствует гиперграни $x_1+\ldots + x_n= -1$. Выберем такое $\alpha>0$, что $c:=\left|\alpha(q_{12}+q_{13}+\ldots+q_{1,2^n-r})\right|<1$ и пусть $d:=q_{1,2^n-r+1}+q_{1,2^n-r+2}+\ldots+q_{1,2^n}$. По сделанному нами выбору в последних $r$ столбцах чисел -1 больше, чем 1, поэтому в виду симметричности неравенств, задающих октаэдр, $d\leq -1$. Положим $\textbf{t}:=(t_1, t_2,\ldots,t_{2^n-r}, t_{2^n-r+1}, \ldots, t_{2^n})=(-c-d,\alpha,\ldots,\alpha,1,\ldots,1)$. Ввиду всевозможных симметрий $\textbf{t}$ является решением не только первого уравнения, но и всей системы.
\end{proof}

Вопрос о реализации других правильных многогранников (их остаётся пять) является сугубо вычислительным, и поэтому мы его опускаем.

\medskip

Напоследок в этом разделе сделаем ещё наблюдения, связывающие свойство ортантности полиэдров разных размерностей.
\begin{proposition}\label{ph1}Пусть полиэдр $P$ с $m$ гипергранями --- ортантный, а $M$ --- такая его гипергрань, что содержит $m-1$ граней полиэдра $P$ коразмерности 2. Тогда $M$ --- ортантный полиэдр.\end{proposition}
\begin{proof}В этих условиях реализация $M$ как сечения ортанта $\R^{m-1}_{\geq 0}$ индуцирована реализацией $P$ как сечения ортанта $\R^m_{\geq 0}$. \end{proof}

\begin{proposition}\label{ph2}
Пусть $\y\in\Y_n(m)$ --- ёжик с иголками $\textbf{a}_1,\ldots,\textbf{a}_m$ такой, что иголки $\textbf{a}_1,\ldots,\textbf{a}_{m-1}$ лежит в одной гиперплоскости $\alpha$. В этом случае эти иголки составляют ёжик $\y'\in\Y_{n-1}(m-1)$. Ёжик $\y$ является ортантным, если и только если $\textbf{a}_m\bot \alpha$ и $\y'$ ортантный.\end{proposition}
\begin{proof}
Можно считать, что $\alpha$ задаётся уравнением $x_n=0$, то есть все последние координаты векторов $\textbf{a}_i$ нулевые при $i\leq m-1$. Запишем $\textbf{a}_m=(a_{m1},\,\ldots,\,a_{mn})$, причём ясно, что $a_{mn}\ne 0$. Рассмотрим соответствующую систему $(!)$. Вспомним, что уравнения этой системы естественно нумеруются парами $(i,j)$, где $1\leq i\leq j\leq n$. В рассматриваемой ситуации уравнения с номерами $(i,n)$, выглядят так: $a_{mi}a_{mn}t_m=\delta_{in}$, где $\delta_{in}$ --- символ Кронекера. Поэтому система $(!)$ имеет положительное решение тогда и только тогда, когда все $a_{mi}=0$, и система, полученная из (!) отбрасыванием уравнений с номерами $(i,n)$ имеет положительное решение, что в свою очередь эквивалентно тому, что $\textbf{a}_m\bot \alpha$, а $\y'$ ортантный.
\end{proof}

\section{Классификация ортантных полиэдров размерности n и ранга n+1}\label{rs}

В этом разделе мы даём полную классификацию ортантных полиэдров в $\Pol^{n+1}_n$. Согласно результатам предыдущего раздела достаточно классифицировать все основные ортантные полиэдры ранга $n+1$, то есть найти критерий ортантности полиэдров из $\Pol^{n+1}_n(n+1)$.

Пусть $P\in \Pol^{n+1}_n(n+1)$ задан системой $(\ast)$. Для таких полиэдров есть ровно две взаимоисключающие возможности:
\begin{itemize}
\item[(1)] Среди векторов $\textbf{a}_i$ найдутся $n$ линейно зависимых;
\item[(2)] Найдётся симплекс $\Delta\sim P$.
\end{itemize}
Согласно предложению~\ref{ph2} первый случай сводится к меньшей размерности, и поэтому достаточно провести классификацию ортантных симплексов. Для того, чтобы выполнить эту классификацию, напомним о понятии ортоцентрических симплексов.

\begin{definition}Симплекс называется \emph{ортоцентрическим}, если его высоты пересекаются в одной точке.\end{definition}

Следующее утверждение хорошо известно, и его доказательство не представляет труда.
\begin{proposition}\label{psimplex}
Следующие условия эквивалентны:
\begin{itemize}
\item[(i)] Симплекс $\Delta$ является ортоцентрическим;
\item[(ii)] Рёбра симплекса $\Delta$ перпендикулярны противолежащим граням коразмерности 2;
\item[(iii)] Все трёхмерные грани симплекса $\Delta$ являются ортоцентрическими тетраэдрами;
\item[(iv)] Для любых четырёх вершин $A$, $B$, $C$, $D$ симплекса $\Delta$ имеет место равенство $AB^2+CD^2=AD^2+BC^2=AC^2+BD^2$.
\end{itemize}
\end{proposition}

\begin{remark}
Ортоцентрический симплекс является в разных смыслах более правильным обобщением понятия треугольника, чем просто симплекс. В частности, как и треугольник, ортоцентрический симплекс может быть \emph{остроугольным}, \emph{прямоугольным} или \emph{тупоугольным}. В ортоцентрическом симплексе при всех вершинах кроме, быть может, одной все двумерные грани имеют острый угол. При оставшейся вершине все двумерные грани одновременно либо имеют острый угол (и тогда симплекс называется остроугольным), либо прямой (прямоугольный симплекс), либо тупой (тупоугольный симплекс). В произвольном же симплексе может оказаться так, что при одной и той же вершине часть двумерных граней имеет острый угол, часть - прямой, часть - тупой.
\end{remark}

\begin{theorem}\label{ts}
Симплекс размерности $n$ может быть реализован как сечение неотрицательного ортанта $\R^{n+1}$ гиперплоскостью тогда и только тогда, когда является остроугольным ортоцентрическим симплексом.
\end{theorem}
\begin{proof} $\underline{\Rightarrow}$ Для симплекса, высеченного гиперплоскостью в $\R^{n+1}$, выполнение пунктов (ii) и (iv) предложения~\ref{psimplex} совершенно очевидно. Остроугольность двумерных граней следует из многократного применения предложения~\ref{ph1}.

$\underline{\Leftarrow}$ Пусть $\frac{1}{\alpha_1}x_1+\ldots+\frac{1}{\alpha_{n+1}}x_{n+1}=1$ --- гиперплоскость в $\R^{n+1}$, где $\alpha_i>0$ для всех $i$. Она высекает $n$-мерный симплекс в $\R^{n+1}_{\geq 0}$, вершины которого имеют координаты $(\alpha_1,0,\ldots,0),\ldots,(0,\ldots,0,\alpha_{n+1})$. Длина рёбра, соединяющего вершины с номерами $i$ и $j$ равна $\sqrt{\alpha_i^2+\alpha_j^2}$.

Пусть теперь $\Delta=A_1\ldots A_{n+1}$ --- остроугольный ортоцентрический симплекс. Тогда рёбра $A_iA_j$ этого симплекса удовлетворяют соотношениям из пункта (iv) предложения~\ref{psimplex} и $A_iA_j^2+A_jA_k^2>A_iA_k^2$. Так как симплекс определяется длинами своих рёбер, ортантность $\Delta$ мы проверим, если докажем, что в этих условиях система из $C^2_{n+1}$ уравнений $x_i+x_j=A_iA_j^2$ имеет положительное решение. А это решение можно просто явно выписать: $x_i=\frac{1}{2}(A_iA_j^2+A_iA_k^2-A_jA_k^2)$.
\end{proof}

\section{Реализация полиэдров сечением ортанта произвольной размерности}\label{rqb}
В этом разделе мы даём (почти полный) ответ на вопрос~\ref{mainq}.B из введения. Полиэдры, для которых этот ответ положительный, в этом разделе мы называем \emph{реализуемыми}.

\begin{theorem}\label{tqb}
 Пусть для некоторого $s\in \N$ конус рецессии $C(P)$ полиэдра $P\in \Pol_n(m)$ вкладывается в $\R^s_{>0}\cup \textbf{0}$. Тогда $P$ реализуемый.
\end{theorem}

Как замечалось в разделе~\ref{intro}, это утверждение почти что даёт критерий о возможности реализации полиэдров указанным способом, ибо, если полиэдр реализуется как сечение плоскостью ортанта $\R^k$, то его конус рецессии автоматически вкладывается в $\R^k_{\geq 0}$, и поэтому открытым вопрос~\ref{mainq}.B остаётся лишь для тех полиэдров, чей конус рецессии вкладывается в $\R^s_{>0}\cup \textbf{0}$, но не в $\R^t_{\geq 0}$.

Идею доказательства теоремы~\ref{tqb} проще всего проиллюстрировать на многогранниках. Основным является следующий факт.
\begin{proposition}\label{pqb} Пусть $P,Q\in\Pol$, $P$ --- многогранник, $P\prec Q$. Тогда, если $Q$ --- ортантный, то $P$ --- реализуемый.
\end{proposition}
\begin{proof}Пусть $Q\in\Pol_n(m)$ задан системой $(\ast)$. Согласно замечанию~\ref{ziez} после некоторого поворота каждой гиперграни $P$ можно сопоставить параллельную ей гипергрань $Q$. При этом можно считать, что разным гиперграням $P$ мы сопоставили разные гиперграни $Q$ (полиэдр $Q$ можно подправить с сохранением свойства ортантности, воспользовавшись преобразованием IV системы $(\ast)$, см. раздел~\ref{rez}). Пусть для определённости гиперграням $P$ сопоставлены гиперграни $Q$, соответствующие первым $m'$ неравенствам в $(\ast)$. Изменим в неравенствах системы $(\ast)$ свободные члены так, чтобы она задавала многогранник $P$: первые $m'$ неравенств подправим так, чтобы они честно задавали $P$, а в остальных неравенствах просто подберём свободные члены так, чтобы им удовлетворяли все точки $P$ (это возможно ввиду ограниченности $P$). Соответствующая система $(!)$ при этом не изменится, поэтому отображение
$$(x_1,\ldots,x_n)\mapsto(k_1(a_{11}x_1+\ldots+a_{1n}x_n-b_{1}),\ldots,k_{m}(a_{m1}x_1+\ldots+a_{mn}x_n-b_{m})),$$
где $k_i^2$ --- решение системы $(!)$, будет изометрией. Эта изометрия и задаст $P$ как сечение $\R^m_{\geq 0}$.\end{proof}

\begin{corollary}\label{cqb}Все многогранники являются реализуемыми.\end{corollary}
\begin{proof}Пусть $P\in\Pol_n$ --- многогранник, а $Q\in\Pol_n^{C^{2}_{n+1}}$ --- ортантный полиэдр (например, из примера~\ref{endgo}). Обозначим $\y:=\yz(P)\cup\yz(Q)$ --- ёжик, полученный объединением иголок на любом общем изображении. Тогда $Q\prec s(\y)$ и $s(\y)$ --- ортантный по следствию~\ref{coo}. Кроме того, $P\prec s(\y)$, и остаётся лишь воспользоваться предложением~\ref{pqb}.\end{proof}

Если предложение~\ref{pqb} попытаться обобщить на случай произвольного полиэдра $P$, то проблема возникнет в тот момент, когда мы будем подбирать свободные члены в последних $m-m'$ неравенствах системы $(\ast)$. Тем не менее, если $P$ удовлетворяет условиям теоремы~\ref{tqb}, для доказательства реализуемости $P$ эту проблему удаётся обойти, взяв в доказательстве следствия~\ref{cqb} ортантный полиэдр $Q\in\Pol_n^{C^{2}_{n+1}}$ специального вида, для которого предложение~\ref{pqb} останется верным.

\begin{proof}[Доказательство теоремы~\ref{tqb}]
Если $C(P)\subset \R^s_{>0}\cup \textbf{0}$, то можно выбрать такое $\eps>0$, что линейные функционалы $x_i-\eps x_j$ положительны на $C(P)$ для любых $i\ne j$. При надлежащем выборе столбца свободных членов $\textbf{b}$, система неравенств $x_i\geq b_i$, $x_i-\eps x_j \geq b^-_{ij}$, $x_i-\eps x_j \geq b^+_{ij}$ доставляет минимальную систему неравенств $(\ast_Q)$, задающих некоторый полиэдр $Q$ в $\R^s$. Несложно убедиться, что $Q$ --- ортантный ранга $C^{2}_{s+1}$. Отметим, что столбец свободных членов может быть изменён так, что точки полиэдра $P$ (при соответствующем вложении в $\R^s$) будут удовлетворять всем неравенствам этой системы. Остаётся добавить к системе $(\ast_Q)$ неравенства, задающие $P$ в $\R^s$, опять заменив столбец свободных членов так, чтобы новая система была минимальной системой неравенств, задающих некоторый (ортантный по следствию~\ref{coo}) полиэдр. Дальнейшая аргументация аналогична предложению~\ref{pqb} и следствию~\ref{cqb}. Отметим, что нас тут не смущает то, что размерность $P$ может оказаться меньше размерности~$Q$.
\end{proof}

\section{Вопросы без ответов}\label{rq}

В этом разделе мы приводим те из вопросов, связанных с проведённым исследованием и оставшихся неразобранными в настоящей работе, что нам кажутся наиболее интересными. Конечно, хотелось бы получить полный ответ на вопрос~\ref{mainq}.A, однако, учитывая, что в настоящей работе в старших размерностях у нас к этому не получилось даже приблизиться, ниже мы приводим более простые вопросы. Получение ответов на вопросы~\ref{q1} и~\ref{q2} стало бы существенным продвижением в сторону решения общей проблемы.

\medskip

В разделе~\ref{rs} мы показали, что метрические условия, которые накладывает условие ортантности на полиэдры из $\Pol_3^4(4)$, имеют простую геометрическую интерпретацию: ортоцентричность симплекса. Интересно, можно ли найти геометрическую интерпретацию в общем случае. Напомним, что метрические соотношения происходят из несколько более слабого условия, чем условие ортантности, а именно из условия, что система (!) имеет решение (необязательно положительное). Иначе говоря, $\rk Q=\rk Q'$, где $Q$ и $Q'$ --- это матрица и, соответственно, расширенная матрица системы (!).
\begin{question}\label{q1}
Какие метрические соотношения накладывает на полиэдр условие $\rk Q=\rk Q'$? Можно ли их задать геометрически?
\end{question}
Особо интересен этот вопрос в случае полиэдров из $\Pol_3^5(5)$, который легко сводится к случаю четырёхугольной пирамиды. Ответ на этот вопрос продвинул бы нас в сторону классификации ортантных полиэдров размерности 3. Напомним, что для завершения этой классификации нам необходимо сделать следующее.
\begin{question}\label{q2}
Описать основные ортантные трёхмерные полиэдры рангов 5 и 6.
\end{question}

\medskip

Напоследок заметим, что деятельность, которую мы развели в этой работе, можно обобщить, задав вопрос, являющийся промежуточным между вопросами~\ref{mainq}.A и~\ref{mainq}.B.

\begin{question}\label{q3}
Пусть $f: \N\to\N$ такая, что $f(m)\geq m$ для всех $m\in\N$, и $P\in\Pol_n(m)$. Можно ли реализовать $P$ как сечение положительного ортанта $\R^{f(m)}_{\geq 0}$?
\end{question}

Ясно, что вопросы~\ref{mainq}.A и~\ref{mainq}.B являются частными случаями вопроса~\ref{q3}. Интересно было бы понять что-то и про другие частные случаи, например, $f(m)=m+1$ или, более общо, $f(m)=m+k$.

Другой способ обобщения нашей деятельности состоит в замене положительного ортанта на другой полиэдральный конус. Например, можно каждому $m\in\N$ поставить в соответствие некоторый полиэдральный конус $\Omega_m\subset\R^m$, и задаться вопросом, какие полиэдры реализуются сечением таких конусов.


\end{document}